\documentclass{amsart}
\usepackage{amssymb}
\usepackage{amsmath}
\usepackage{latexsym}%
\usepackage[normalem]{ulem}
\usepackage{color}
\usepackage[all]{xy}
\usepackage[linktoc=all, pagebackref, hyperindex]{hyperref}%
\hypersetup{colorlinks,  citecolor=blue,  filecolor=blue,  linkcolor=red,  urlcolor=black}

\newcommand{\OS}{\bar{S}}

\newcommand{\s}{\mathbf{s}}

\newcommand{\m}{\mathfrak{m}}

\def\RR{{\mathbb{R}}}
\def\NN{{\mathbb{N}}}
\def\QQ{{\mathbb{Q}}}
\def\ZZ{{\mathbb{Z}}}
\def\KK{{\mathbb{K}}}

\DeclareMathOperator{\Hom}{Hom}
\DeclareMathOperator{\Max}{Max}

\DeclareMathOperator{\gr}{group}

\DeclareMathOperator{\Ap}{Ap}

\newcommand{\AP}{\operatorname{Ap}}

\newcommand{\PF}{\operatorname{PF}}
\newcommand{\QF}{\operatorname{QF}}
\newcommand{\co}{\operatorname{cone}}

\newcommand{\ord}{\operatorname{ord}}
\newcommand{\Ext}{\operatorname{Ext}}
\newcommand{\typ}{\operatorname{type}}
\newcommand{\Gr}{\operatorname{group}}
\newcommand{\mgs}{\operatorname{mgs}}
\newcommand{\depth}{\operatorname{depth}}
\newcommand{\relint}{\operatorname{relint}}
\def\c{{\mathbf c}}
\def\m{{\mathbf m}}
\def\1{{\mathbf 1}}
\def\a{{\mathbf a}}
\def\b{{\mathbf b}}
\def\z{{\mathbf z}}
\def\f{{\mathbf f}}

\def\bs{{\mathbf s}}
\def\x{{\mathbf x}}

\def\v{{\mathbf v}}
\def\w{{\mathbf w}}

\def\fc{{\mathfrak c}}
\def\fm{{\mathfrak m}}

\def\fp{{\mathfrak p}}
\def\fq{{\mathfrak q}}

\newtheorem{theorem}{Theorem}[section]
\newtheorem{corollary}[theorem]{Corollary}
\newtheorem{lemma}[theorem]{Lemma}
\newtheorem{proposition}[theorem]{Proposition}

\theoremstyle{definition}
\newtheorem{definition}[theorem]{Definition}
 
\newtheorem{example}[theorem]{Example}
\newtheorem{remark}[theorem]{Remark}

\numberwithin{equation}{section}

\title{Type and conductor of simplicial affine semigroups}

\author{Raheleh Jafari}
\address{Mosaheb Institute of Mathematics, Kharazmi University, 
	Tehran, Iran}
\email{rjafari@ipm.ir}

\author{Marjan Yaghmaei}
\address{Faculty of Mathematical Sciences and Computer, Kharazmi University, Tehran, Iran}
\email{hasti.tmu83@yahoo.com}

\begin{document}

\begin{abstract}
	  We provide a generalization of  pseudo-Frobenius numbers of numerical semigroups to the context  of the simplicial affine semigroups. In this way, we characterize the Cohen-Macaulay type of the simplicial affine semigroup ring $\KK[S]$.
 We define the type of $S$, $\typ(S)$,  in terms of some Ap\'ery sets of $S$ and show that it coincides with the Cohen-Macaulay type of the semigroup ring, when $\KK[S]$ is Cohen-Macaulay.  If  $\KK[S]$ is a  $d$-dimensional Cohen-Macaulay ring of embedding dimension at most $d+2$, then $\typ(S)\leq 2$. Otherwise,   $\typ(S)$ might be arbitrary large and it has no upper bound in terms of the embedding dimension.  Finally, we present a generating set for the conductor  of $S$ as an ideal of its normalization. 
 	\end{abstract}


\maketitle

	\section{Introduction}
	Let $S$ be an affine semigroup in $\NN^d$, where $\NN$ denotes the set of nonnegative integers. The 
affine semigroup  ring $\KK[S]$, over a field $\KK$, is defined as the subring 
$\{\oplus_{\a\in S}k_\a\x^\a \colon k_\a\in\KK\}$
	of the polynomial  ring $\KK[\x]=\KK[x_1,\dots,x_d]$.
		The structure of  $\KK[S]$ is
	intimately related to the structure of the affine semigroup $S$ and $\co(S)$, the rational polyhedral cone spanned by $S$.
The Cohen-Macaulay,  Gorenstein and Buchsbaum properties of $\KK[S]$ have been characterized in terms of certain numerical and topological properties of $S$, see  \cite{Goto-Watanabe-1978}, \cite{Hoa-Trung-1986}, \cite{Hoa-Trung-1988}, \cite{GSW-76}, \cite{RG-1998} and \cite{GR-2002}.  Our aim in this paper, is to characterize the conductor and Cohen-Macaulay type of $\KK[S]$.

If $d=1$, then $S$ is a submonoid of $\NN$. Let $h$ be the greatest common divisor of nonzero elements in $S$. Dividing all elements of $S$ by $h$, we obtain an isomorphic  semigroup in $\NN$. 
A submonoid $S$ of $\NN$ such that $\gcd(s ; s\in S)=1$ is called a {\em numerical semigroup}.  In other words, the study of affine semigroups in $\NN$ reduces to the study of numerical semigroups.  The condition $\gcd(s ; s\in S)=1$ is equivalent to  say that $\NN\setminus S$ is a finite set, \cite[Lemma~2.1]{RG}.
 For an affine semigroup $S$, consider 
the natural partial ordering $\preceq_S$ on $\NN^d$ where, for all elements $x,y\in\NN^d$, $x\preceq_Sy$ if $y-x\in S$.  For  a numerical semigroup $S\subsetneqq\NN$, the maximal elements of $\NN\setminus S$ with respect to $\preceq_S$ are called {\em pseudo-Frobenius numbers}. 
	Fr\"{o}berg, Gottlieb and H\"{a}ggkvist \cite{FGH}, defined the type of the numerical semigroup $S$ as the
	cardinality of the set of its pseudo-Frobenius numbers. This notion of type coincides with the Cohen-Macaulay type of the  numerical semigroup ring $\KK[S]$, see \cite{Stamate-2018} for a detailed proof. 
	
By analogy, Garc\'{\i}a-Garc\'{\i}a, Ojeda, Rosales and  Vingneron-Tenorio, define a pseudo-Frobenius element of $S$ to be an element $\a\in\NN^d\setminus S$ such that $\a+S\setminus\{0\}\subseteq S$, in \cite{GORV}. They show that the set of pseudo-Frobenius elements of $S$, $\PF(S)$, is not empty, precisely when $\depth\KK[S]=1$. Thus, when $d>1$ and $\KK[S]$ is a Cohen-Macaulay ring, the set of pseudo-Frobenius elements of $S$ is empty and express noting about the Cohen-Macaulay type of the semigroup ring.

In Section~\ref{sec:type}, we present another generalization of pseudo-Frobenius numbers, determining  the Cohen-Macaulay type of the semigroup ring $\KK[S]$, under the assumption that the affine semigroup $S\subset \NN^d$ is also simplicial, i.e. $\co(S)$ has $d$ extremal rays.  All  fully embedded affine semigroups in $\NN^d$, for $d=1,2$, are simplicial. Let $\a_1,\dots,\a_d$ be the componentwise smallest integer vectors in $S$, situated on each extremal ray of $\co(S)$, respectively. Then $\cap^d_{i=1}\Ap(S,\a_i)$ is a finite set, where $\Ap(S,\a_i)=\{\a\in S ; \a-\a_i\notin S\}$ denotes the Ap\'ery set of $S$ with respect to $\a_i$. In numerical case, where $d=1$, $\a_1$ is the smallest positive number in $S$ and the set of pseudo-Frobenius numbers is equal to the set $\{\w-\a_1 ; \w\in\Max_{\preceq_S}\Ap(S,\a_1)\}$. In an analogue way, we consider the set 
$$\QF(S)=\{\w-\sum^d_{i=1}\a_i \ ; \ \w\in\max_{\preceq_S}\cap^d_{i=1}\Ap(S,\a_i)\}.$$
We call the elements of $\QF(S)$, {\em quasi-Frobenius} elements, to distinguish from $\PF(S)$ introduced in \cite{GORV}. 
In Proposition~\ref{type}, we show that the cardinality of $\QF(S)$ which is called the {\em type} of $S$ and is denoted by $\typ(S)$, equals the Cohen-Macaulay type of $\KK[S]$, when the semigroup ring is Cohen-Macaulay. Campillo and Gimenez in \cite{Campillo-Gimenez-2000}, study combinatorics of some simplicial complexes associated to elements of $S$. They compute the homologies of the simplicial complexes in terms of certain graph homologies, and show in \cite[Theorem~4.2(ii)]{Campillo-Gimenez-2000}, that for Cohen-Macaulay simplicial affine semigroups, $\typ(S)$ equals the  number of maximal elements of $\Ap(S,E)$ with respect to $\preceq_S$. Our result, Proposition~\ref{type}, provides a different algebraic proof for this fact. 
 Although $\PF(S)$ and $\QF(S)$ coincide  in the numerical case where $d=1$, but for $d>1$ they have no common element, see Remark~\ref{rem:QFPF}.  

As an immediate consequence, we derive that the simplicial affine semigroup ring $\KK[S]$ is Gorenstein if and only if  it is Cohen-Macaulay and $\cap^d_{i=1}\Ap(S,\a_i)$ has a single maximal element with respect to $\preceq_S$. This result is already proved in \cite{Campillo-Gimenez-2000,RG-1998}, by different  arguments. In the rest of Section~\ref{sec:type}, we investigate the above bound for $\typ(S)$. Generalizing \cite[Theorem~11]{FGH}, we show in Theorem~\ref{aa}, that $\typ(S)\leq2$, if $S$ is generated by $d+2$ elements and $\KK[S]$ is Cohen-Macaulay. If either $\KK[S]$ is not Cohen-Macaulay or if the minimal generating set of $S$ has more than $d+2$ elements, then $\typ(S)$ might be arbitrary large, see Example~\ref{3^n} and Example~\ref{ex:3.7}.

 Recall that the normalization of an integral domain $R$ is the set of elements
	in its field  of fractions satisfying a monic polynomial in $R[y]$. Then 
  $R=\KK[S]$ is an integral domain with  normalization  $\bar{R}=\KK[\Gr(S)\cap\co(S)]$, where $\Gr(S)$ denotes the group of differences of $S$, \cite[Proposition~7.25]{Miller-Sturmfels}. The purpose of the last section, is  to  investigate the normality  of  $R$ and, to  detect a generating set for the  conductor of $R$, $C_R=(R:_T\bar{R})$,  where $T$ denotes the total ring of fractions of $R$.   We will work with corresponding objects in $S$. In semigroup interpretation,  $\bar{S}=\Gr(S)\cap\co(S)$ and $\fc(S)=\{\b\in S \ ; \ \b+\bar{S}\subseteq S\}$, are called the   {\em normalization} and the  {\em conductor} of $S$, respectively, \cite{Bruns-Gubeladze-2009}. The ring $R$ is normal precisely when the semigroup  $S$ is normal, i.e. $S=\bar{S}$, \cite{Bruns-Herzog}. 
 
 Quasi-Frobenius elements, are also profitable to recognize the normality of $S$. 	 Note that quasi-Frobenius elements, might have negative components. Having more negative components in elements of $\QF(S)$, makes  the semigroup more  close  to being normal. More precisely, if $S$ is normal, then $-\QF(S)\subseteq\co(S)$ and so $-\QF(S)\subseteq\NN^d$. The converse holds, if $-\QF(S)$ is a subset of the relative interior of $\co(S)$. This is the subject of Theorem~\ref{QF&NORMAL}, which states also that in this case $-\QF(S)\subseteq S$.
 	
 If $S\subseteq\NN$ is a  numerical semigroup, then $\bar{S}=\NN$ and $\fc(S)$ is a principal ideal of $\NN$ generated by $F+1$, where $F=\max(\ZZ\setminus S)$. The rest of Section~\ref{sec:norm}, is devoted to find a generating set for $\fc(S)$, as an ideal of $\bar{S}$, where $S$ is an arbitrary simplicial affine semigroup. 
 
  Note that any element $\c\in\co(S)\cap\NN^d$, is uniquely presented as $\c=\sum^d_{i=1}n_i\a_i+r(\c)$  for some $n_1,\dots,n_d\in\NN$ and $r(\c)\in P_S$,  where $P_S=\{\sum^d_{i=1}r_i\a_i \ ; \ 0\leq r_1,\dots,r_d<1\}$  is the fundamental parallelogram of $S$. It is not difficult to observe that  $P_S\cap\Gr(S)=\{ r(\w) \ ; \ \w\in \cap^d_{i=1}\Ap(S,\a_i)\}$, see Lemma~\ref{lem:OS}. Let $P_S\cap\Gr(S)=\{0=\b_0,\b_1,\dots,\b_k\}$.
 For $\w_1,\dots,\w_k\in \cap^d_{i=1}\Ap(S,\a_i)$ with $r(\w_i)=\b_i$, we consider the vector $\f_{(\w_1,\dots,\w_k)}=\sum^d_{i=1}f_i\a_i$, where $f_i$  is the maximum integer such that $\w_j-f_i\a_i\in P_S$ for  some $j=1,\dots,k$. 
  We show in Theorem~\ref{fff}, that any element of the minimal generating set for $\fc(S)$, as an ideal of $\bar{S}$, is of the form 
 $\f_{(\w_{1},\dots,\w_{k})}-\b_j+\sum^d_{i=1}l_i\a_i$ for some $l_i\in\{0,1\}$ ,  $j\in\{1,\ldots,k\}$    
  and
  $\w_1,\dots,\w_k\in\cap^d_{i=1}\Ap(S,\a_i)$ such that $r(\w_i)=\b_i$ for $i=1,\dots,k$. Moreover, at least for one $i$, we have $l_i=0$. As $\cap^d_{i=1}\Ap(S,\a_i)$ is a finite set and $\fc(S)=\{\c\in S \ ; \ \c+\b_i\in S \text{ for } i=1,\dots,k\}$,  Theorem~\ref{fff} provides an algorithmic way to find a generating set for $\fc(S)$. 
Let $\preceq_c$ denote the natural partial ordering with respect to $\co(S)$, that is $\a\preceq_c\b$ if $\b-\a\in\co(S)$. 
We show in Corollary~\ref{cor:principal} that, if $\KK[S]$ is Cohen-Macaulay  and   $\max_{\preceq_{c}}\cap^d_{i=1}\Ap(S,\a_i)=\{\w\}$, then  	$\{\w-\b \ ; \b\in\max_{\preceq_c}r(\Ap(S,E)) \},$ generates $\fc(S)$ as an ideal of $\bar{S}$.
  Several explicit examples are provided to illustrate the generating set of $\fc(S)$, as an ideal of $\bar{S}$.

	\section{Fundamentals}\label{sec:Back}
	\label{background}
 By an affine semigroup, we mean a finitely generated submonoid of $\NN^d$, where $\NN$ denotes the set of nonnegative integers and $d\in\NN\setminus\{0\}$. Let $S$ be an affine semigroup minimally generated by $\mgs(S)=\{\a_1,\dots,\a_e\}$. We  write $S=\langle\a_1,\dots,\a_e\rangle$, to indicate its generating set. The minimal generating set, $\mgs(S)$, is a unique finite set, see \cite[Chapter~3]{RG-1999}. The number of elements in $\mgs(S)$ is called the {\em embedding dimension} of $S$.  For a field  $\KK$, the semigroup ring $\KK[S]$ is the subalgebra of the polynomial ring $\KK[x_1,\dots,x_d]$ generated by the monomials with exponents in $S$. The ring $\KK[S]=\KK[\x^{\a_1},\dots,\x^{\a_e}]$ has a unique maximal monomial ideal $\fm=(\x^{\a_1},\dots,\x^{\a_e})$.

 Given two sets $A,B\subseteq\NN^d$, we write $A+B$ for the set $\{\a+\b \ ; \ \a\in A , \b\in B\}$. If $A=\{\a\}$, we simply write $\a+B$, instead of $\{\a\}+B$. 
Recall that an {\em ideal} of $S$ is a nonempty set $H\subseteq S$, such that $H+S\subseteq S$.
 For any ideal $H$ of $S$, there exists a set of vectors $B=\{\b_1,\dots,\b_l\}$ such that $H=B+S=\cup^l_{i=1}\b_i+S$. In this case, $B$ is called a {\em generating set} of $H$. If no proper subset of $B$ generates $H$, we refer to $B$ as the minimal generating set of $H$.  Any ideal of an affine semigroup has a unique minimal generating set.
 Note that for ideals $H_1$ and $H_2$ of $S$, $H_1+H_2$ is also an ideal of $S$. In particular,  for an ideal $H$ of $S$, $n$ times summation of $H$ which is denoted by   $nH$,  is again an ideal of $S$. Let $M=S\setminus\{0\}$ be the maximal ideal of $S$. 

  A  monomial in the semigroup ring $\KK[S]$ is an element of
 	the form $\x^\a$ for $\a\in S$. An ideal $I\subseteq\KK[S]$ is a monomial ideal if it is
 	generated by monomials. For any subset $H$ of $S$, let  $\KK[H]$ denote  the $\KK$-linear span of the
 	monomials $\x^\a$ with $\a\in H$. Then  $I$ is a monomial ideal if and only if $I=\KK[H]$ for some ideal $H$ of $S$, or equivalently, if $I$ is homogeneous with respect
 	to the tautological grading on $\KK[S]$,  which is defined by $\deg(\x^\a)=\a$.  Note that  $\fm=\KK[M]$.

The affine semigroup $S\subseteq\NN^d$ is called  {\em simplicial} if there exist $d$  elements  $\a_{i_1},\dots,\a_{i_d}\in\mgs(S)$  such that they are linearly independent over the field of rational numbers $\QQ$ (equivalently, over the field of real numbers $\RR$), and for each element $\a\in S$, we have $n\a\in\NN\a_{i_1}+\dots+\NN\a_{i_d}$, for some positive integer $n$. There is a geometrical interpretation for the simplicial property. Let 
$$
\co(S)=\left \{ \sum_{i=1}^e \lambda_i \a_i \ ; \ \lambda_i \in \RR_{\geq 0}, \text{ for } i=1,\dots, e \right\},
$$
denote the rational polyhedral cone generated by $S$.
The dimension (or rank) of $S$ is defined as the dimension of  the affine subspace it generates, which is the same as the dimension of the subspace generated by  $\co(S)$.  
 The $\co(S)$ is  polyhedral i.e. it is the intersection of finitely many closed linear
half-spaces in $\RR^d$, each of whose bounding hyperplanes contains the origin, \cite[Corollary~7.1(a)]{Schrijver}. These half-spaces are called {\em support hyperplanes}. The integral vectors in each  support hyperplane, is a face of $S$, and all maximal faces (called facets) are in this form.  The intersection of any two adjacent support hyperplanes is a one-dimensional vector space, which is called an {\em extremal ray}.    The $\co(S)$ has at least $d$ facets and at least $d$ extremal rays. It has $d$ facets (equivalently,  it has $d$ extremal rays), precisely when  $S$ is simplicial.

Throughout this paper, we  assume that $S$ is a simplicial affine semigroup. 
On each extremal ray of $\co(S)$,  the componentwise smallest element from $S$, is  called an {\em  extremal ray} of $S$. 
Denote by $\a_1,\dots, \a_d$ the extremal rays of $S$. These form a basis for $\co(S)$. 
For $\z\in \RR^d$ such that $\z=\sum_{i=1}^d \lambda_i \a_i$ with $\lambda_i \in \QQ, i=1,\dots,d$, we set $[\z]_i=\lambda_i$ for $i=1,\dots, d$.  For each element $\a\in S$, we have $n\a\in\NN\a_{1}+\dots+\NN\a_{d}$, for some positive integer $n$. In other words, $\{\x^{\a_1},\dots,\x^{\a_d}\}$ provides a monomial system of parameters for $\KK[S]$. The {\em fundamental (semi-open) parallelotope} of $S$ is the set
\begin{eqnarray*}
	P_S&=&\left\{ \z \in\RR^d \ ; \  0\leq [\z]_i <1 \text{ for } i=1, \dots, d\right\}\\
	&=&{}\{\sum^d_{i=1}\lambda_i\a_i \ ; \ \lambda_i\in\QQ \ , \ 0\leq\lambda_i<1 \text{ for } i=1, \dots, d\}. 
\end{eqnarray*}
	Its closure in $\RR^d$ is the set $\bar{P}_S=\left\{ \z \in\RR^d \ ; \    0\leq [\z]_i \leq 1 \text{ for } i=1, \dots, d\right\}$.
	It is well known, and easy to see, that any $\a$ in $\co(S)\cap\NN^d$ decomposes uniquely as $\a=\sum_{i=1}^d n_i \a_i+ r(\a)$, with $r(\a) \in P_S \cap \NN^d$ and nonnegative integers $n_1,\dots, n_d$.  We will call $r(\a)$ the  {\it remainder} of $\a$ in $P_S$. 
\begin{remark}\label{e(a)}
	Let $\a\in\co(S)\cap\NN^d$ and let  $n_i=\lfloor[\a]_i\rfloor$ be the unique integer such that $n_i\leq[\a]_i<n_i+1$, for $i=1,\dots,d$.  Then $r(\a)=\sum^d_{i=1}([\a]_i-n_i)\a_i$. 
\end{remark}	

 We consider the natural partial orderings $\preceq_S$  and $\preceq_c$ on $\NN^d$ where, for all elements $\a$ and $\b$ in $\NN^d$, $\b\preceq_S\a$ ($\b\preceq_c\a$), if there is an element $\c\in S$ ($\c\in\co(S)$) such that $\a=\b+\c$.  The partial order  $\preceq_c$ is indeed the  coordinatewise order  on $\co(S)\cap\NN^d$. More precisely, for $\a,\b\in\co(S)\cap\NN^d$,  $\a\preceq_c\b$ if and only if  $[\a]_i\leq[\b]_i$ for $i=1,\dots,d$.

An element $\a\in S$,  may be written as $\a=\sum^{d+r}_{i=1}l_i\a_i$ for some nonnegative integers $l_1,\dots,l_{d+r}$. The value $\sum^{d+r}_{i=1}l_i$ is called the {\em length} of the expression $\sum^{d+r}_{i=1}l_i\a_i$. The maximum integer $n$ such that $nM$ contains $\a$, is called the {\em order} of $\a$ and it is denoted by $\ord(\a)$. In other words, $\a\in nM\setminus (n+1)M$ if and only if $n=\ord(\a)$.  The expression of length $\ord(\a)$ of $\a$,    is called a {\em maximal expression} of $\a$.

The {\em Ap\'{e}ry set} of an element $\b\in S$ is defined as $\Ap(S,\b)=\{\a\in S \ ; \ \a-\b\notin S\}$. We will denote the zero vector of $\NN^d$ by $0$. Since $S\subseteq\NN^d$, for $\b\neq 0$ we have  $0\in \Ap(S,\b)$.    Note that if $\a\in\Ap(S,\b)$ and $\z\in S$ such that $\z\preceq_S\a$, then $\z\in\Ap(S,\b)$. For a subset $E$, we set
 $$\Ap(S,E)=\{\a\in S \ ; \ \a-\b\notin S, \text{ for all } \b\in E\}.$$  Throughout the paper, $E=\{\a_1,\dots,\a_d\}$ will denote the set of extremal rays of $S$. Then $\Ap(S,E)=\cap^d_{i=1}\Ap(S,\a_i)$.

 Let $I_S$ denote the kernel of the $\KK$-algebra homomorphism $\varphi:\KK[z_1,\dots,z_{d+r}]\longrightarrow\KK[S]$, defined by $z_i\mapsto\x^{\a_i}$, for $i=1,\dots,d+r$.  Then $I_S$ is a binomial prime ideal, \cite[Proposition~1.4]{Herzog-1970}. Note that  $\KK[S]\cong\KK[z_1,\dots,z_{d+r}]/I_S$ has  $S$-graded structure defined by $\deg_S(z_1^{n_1}\dots z_{d+r}^{n_{d+r}})=\sum^{d+r}_{i=1}n_i\a_i$. 
  Let $R'=\frac{\KK[z_1,\dots,z_{d+r}]}{I_S+(z_1,\dots,z_d)}$. Then, as a $\KK$-vector space, $R'$ is  generated by the set of monomials $\z^\a$ such that $\a\notin I_S+(z_1,\dots,z_d)$. Let  $B$  denote  the monomial $\KK$-basis of $R'$.  From \cite[Theorem~3.3]{Ojeda-Tenorio}, we have
\[
\Ap(S,E)=\{\deg_S(u) \ ; \ u\in B\}.
\]
Therefore, as an algorithm to find $\Ap(S,E)$, one may first compute $I_S$, using any of computer algebra systems GAP \cite{GAP}, Singular \cite{Singular}, CoCoA \cite{cocoa}  or Macaulay2 \cite{Mac2}, and then, find the  monomial basis  of $\frac{\KK[z_1,\dots,z_{d+r}]}{I_S+(z_1,\dots,z_d)}$.

\begin{example}\label{2.2}
	Let $\a_1=(5,3,1), \a_2=(1,5,2), \a_3=(8,3,5), \a_4=(2,1,1),\a_5=(2,2,1)$. A computation by Macaulay2 \cite{Mac2}, shows that $$I_S=(z_5^5-z_1z_2z_4^2, z_4^{19}-z_1^2z_3^3z_5^2, z_4^{17}z_5^3-z_1^3z_2z_3^3).$$ Consequently,  $I_S+(z_1,z_2,z_3)=(z_1,z_2,z_3,z_4^{19},z_5^{5},z_4^{17}z_5^3)$. The image of 
	\[\{1, z_4^{r}z_5, z_4^{r}z_5^2, z_4^{s}z_5^3, z_4^{s}z_5^4,  z_4^{t} \ ; \ 0\leq r\leq 18, 0\leq s\leq 16 ,1\leq t\leq18 \}
	\] in $\frac{\KK[z_1,\dots,z_5]}{I_S+(z_1,z_2,z_3)}$, provides a $\KK$-basis. Therefore,  $\Ap(S,E)$ is equal to the set
	$$\{0,r\a_4+\a_5, r\a_4+2\a_5, s\a_4+3\a_5,    s\a_4+4\a_5,  t\a_4 \ ; \ 0\leq r\leq 18, 0\leq s\leq 16 ,1\leq t\leq18  \}.$$
\end{example}

  We write $\Gr(S)$ for the {\em group of differences of $S$}, i.e. $\Gr(S)$ is the smallest group (up to isomorphism) that contains $S$.
	\[
	\Gr(S)=\{\a-\b \mid \a, \b \in S\}.
	\]
By $\Gr(\a_1,\dots,\a_d)$, we mean the 	 smallest group that contains $\{\a_1,\dots,\a_d\}$, equivalently  $\Gr(\a_1,\dots,\a_d)=\{\sum^d_{i=1}z_i\a_i \ ; \ z_i\in\ZZ\}$. 
\begin{remark}\label{rmk:GrAp}
For  $\c\in\Gr(S)$, there exists $\b\in S$ such that $\c+\b\in S$. As $S$ is simplicial, $n\b\in\sum^d_{i=1}\NN\a_i$ for some positive integer $n$.    Therefore, $\c+n\b=\c+\sum^d_{i=1}r_i\a_i\in S$, for some $r_i\in\NN$. Consider $r_1,\dots,r_d$, as small as possible with this property, i.e. 
	\begin{equation}\label{min}
	\c+(r_j-1)\a_j+\sum^d_{i=1, i\neq j }r_i\a_i\notin S,
	\end{equation}
	for $j=1,\dots,d$.   Let  $\c+\sum^d_{i=1}r_i\a_i=\w+\sum^d_{i=1}s_i\a_i$ for some $\w \in\Ap(S,E)$ and $s_1,\dots,s_d\in\NN$. If $r_j>0$, for some $1\leq j\leq d$, then 
		\[
		\c+(r_j-1)\a_j+\sum^d_{i=1, i\neq j }r_i\a_i=\w+(s_j-1)\a_j+\sum^d_{i=1, i\neq j }r_i\a_i.
		\]
By (\ref{min}), we get $s_j=0$.  	
\end{remark}

As an advantage of  considering $\Ap(S,E)$, we recall the following criteria for the Cohen-Macaulay property of $\KK[S]$. 

\begin{proposition}\cite[Corollary~1.6]{RG-1998}\label{CM}
	The following statements are equivalent.
	\begin{enumerate}
		\item $\KK[S]$ is Cohen-Macaulay.
		\item For all $\w_1,\w_2\in\Ap(S,E)$, if $\w_1-\w_2\in\Gr(\a_1,\dots,\a_d)$, then $\w_1=\w_2$.
\end{enumerate}	
\end{proposition}

The Cohen-Macaulay property is indeed equivalent to have a one to one correspondence between elements in $\Ap(S,E)$ and their remainders. More precisely, let 
$$r(\Ap(S,E))=\{r(\w) \ ; \ \w\in\Ap(S,E)\}=\{0=\b_0,\b_1,\dots,\b_k\}.$$
 Then the family of subsets $C_j=\{\w\in\Ap(S,E) \ ; \ r(\w)=\b_j\}$, for $ j=0,\dots,k$,  defines a partition of $\AP(S,E)$. 	
	\begin{lemma}\label{CM-Ci} The following statements hold.
	\begin{enumerate}
		\item $\KK[S]$ is Cohen-Macaulay if and only if $C_i$ is a singleton for $i=1,\dots,k$. 
		\item  $\KK[S]$ is  Buchsbaum if and only if, either  $C_i$ is  a singleton or  $C_i=\{\c+\a_1,\dots,\c+\a_d\}$ for some $\c\in\NN^d$ such that $\c+(S\setminus\{0\})\subset S$, for $i=1,\dots,d$. 
			\end{enumerate}
	\end{lemma}
	
	\begin{proof}
		Let $\v,\w\in\Ap(S,E)$. Then $\v-\w=r(\v)+\sum^d_{i=1}r_i\a_i-r(\w)-\sum^d_{i=1}s_i\a_i$, for  some $r_i,s_i\in\NN$. In particular, $\v-\w=r(\v)-r(\w)+\b$, where $\b=\sum^d_{i=1}(r_i-s_i)\a_i\in\gr(\a_1,\dots,\a_d)$. 
		Therefore, $\v-\w\in\gr(\a_1,\dots,\a_d)$ if and only if $r(\v)=r(\w)$, equivalently $\v,\w\in C_i$, for some $0\leq i\leq k$. 	  Thus, $C_i$s are precisely the equivalence classes under the equivalence relation $\sim$ on $\Ap(S,E)$, where $\w_i\sim \w_j$ if $\w_i-\w_j\in\gr(\a_1,\dots,\a_d)$. 
				Now, the statement (1) follows by Proposition~\ref{CM}, and    the statement (2)  is a consequence of \cite[Theorem~11]{GR-2002}.  
	\end{proof}

	When $d=1$, equivalently  $S$ is isomorphic to a numerical semigroup,   then $\KK[S]$ is a one-dimesional domain and consequently it is Cohen-Macaulay. From another point of view, since all elements of $S$  belong to the real line, one may easily  use Proposition~\ref{CM} or Lemma~\ref{CM-Ci} to check that  $\KK[S]$ is Cohen-Macaulay. The following lemma, extends this property to  simplicial affine semigroups in $\NN^d$.

\begin{lemma}\label{rmk:cm} 
	Let $S$ be of embedding dimension $d+r$.
If the vectors $\a_{d+1},\dots,\a_{d+r}$ belong to the same line passing through the origin of coordinates,
then $\KK[S]$ is Cohen-Macaulay.
\end{lemma}

\begin{proof}
 As $\a_{d+1},\dots,\a_{d+r}$ belong to the same line passing through the origin of coordinates, there exist nonnegative rational numbers $l_i\in\QQ$ such that $\a_{d+i}=l_i\a_{d+1}$, for $i=1,\dots,r$. Let $\w_1,\w_2\in\Ap(S,E)$, then $\w_j=\lambda_j\a_{d+1}=\lambda_j(\sum^d_{i=1}\mu_i\a_i)$, for some nonnegative rational numbers $\lambda_1,\lambda_2,\mu_1,\dots,\mu_d$. 
 If $\w_1-\w_2\in\Gr(\a_1,\dots,\a_d)$, then 
 \[
 (\lambda_1-\lambda_2)(\sum^d_{i=1}\mu_i\a_i)=\w_1-\w_2=\sum^d_{i=1}z_i\a_i,
  \]
 for some integers $z_1,\dots,z_d$. Hence, for all $i$, $z_i=(\lambda_1-\lambda_2)\mu_i$, which forces the signs of $z_i$
 to be the same for all $i$. If they are all nonnegative (nonpositive), the equation $\w_1=\w_2+\sum^d_{i=1}z_i\a_i$ ($\w_2=\w_1+\sum^d_{i=1}(-z_i)\a_i$),  implies $z_i=0$,  since $\w_1,\w_2\in\Ap(S,E)$. 
 Therefore, $\KK[S]$ is Cohen-Macaulay by Proposition~\ref{CM}.
\end{proof}

 The following lemma states an easy but useful property about maximal expressions of Ap\'ery elements in $\Ap(S,E)$, when $S$ has only two nonextremal generators that belong to  the same line passing through the origin of coordinates.

\begin{lemma}\label{edim2}
	Let $S$ be of embedding dimension $d+2$. Then 
	\begin{enumerate}
		\item There are no elements in $\Ap(S,E)$ having two different expressions with the same length.  In particular, each  element in $\AP(S,E)$ has a unique maximal expression.
		\item Assume that $r\a_{d+1}=s\a_{d+2}$, where $r$ and $s$ are relatively prime positive integers with $r>s$. For $\b\in\Ap(S,E)$, an expression $\b=n_1\a_{d+1}+n_2\a_{d+2}$ is maximal if and only if  $n_2<s$.	
	\end{enumerate}
\end{lemma}
\begin{proof}
(1). Let  $\mgs(S)=\{ \a_1,\dots,\a_d,\a_{d+1},\a_{d+2}\}$. 
	Assume on the contrary that, an element   $\b\in\AP(S,E)$ has two  expressions  of the same length $$\b=n_1\a_{d+1}+n_2\a_{d+2}=n_1^{\prime}\a_{d+1}+
	n_2^{\prime}\a_{d+2}.$$
	Then 
	$(n_1-n_1^{\prime})\a_{d+1}+(n_2-n_2^{\prime})\a_{d+2}=0$.
	Since $n_1+n_2=n_1^{\prime}+n_2^{\prime}$,
	$n_1-n_1^{\prime}=n_2^{\prime}-n_2\neq 0$. Therefore $\a_{d+1}=-\a_{d+2}$,  a contradiction.

(2).
 If $n_2\geq s$, then $\b=(n_1+r)\a_{d+1}+(n_2-s)a_{d+2}$. Since $(n_1+r+n_2-s)>n_1+n_2$, the expression $n_1\a_{d+1}+n_2\a_{d+2}$ is not maximal.

Now assume that $n_2<s$.
Let $\b=n_1^{\prime}\a_{d+1}+n_2^{\prime}
	\a_{d+2}$ be a maximal expression. Then $n_2^{\prime}$ is also smaller than $s$, by our  first argument.  Comparing the two expressions of $\b$, we derive
$$(n_1-n_1^{\prime})\a_{d+1}=(n_2^{\prime}-n_2)\a_{d+2}=\frac{r(n_2^{\prime}-n_2)}{s}\a_{d+1}.$$ Since $r$ and $s$ are relatively prime, this follows  $(n_2^{\prime}-n_2)$ is a multiple of $s$, which contradicts $|n'_2-n_2|<s$,  unless $n'_2=n_2$.	
\end{proof}

\section{The type of simplicial affine semigroups}\label{sec:type}

 In this section, we are looking for a characterization of the Cohen-Macaulay type of the affine semigroup ring $\KK[S]$, in terms of some numerical invariants of $S$. All over  the section, as thorough the paper,  
 $S$  is a $d$-dimensional  simplicial affine semigroup with $\mgs(S)=\{\a_1,\dots,\a_{d+r}\}$,  where  $\a_1,\dots,\a_d$ are the  extremal rays of $S$.
 
 If $d=1$, then  dividing elements of $S$ by the greatest common divisor of $\a_1$,$\dots$, $\a_{d+r},$ we obtain an isomorphic semigroup. So, we may assume that $S$ is a numerical semigroup, equivalently $S$ is a submonid of $\NN$ such that $\NN\setminus S$ is a finite set. If $S\neq\NN$, the elements of $\PF(S):=\max_{\preceq_S}\NN\setminus S=\{ f\in\NN\setminus S \ ; \ f+s\in S, \text{ for all } s\in S\setminus\{0\}\}$ are called  pseudo-Frobenius numbers of $S$. The cardinality of $\PF(S)$ is called  the type of $S$, by Fr\"{o}berg, Gottlieb and H\"{a}ggkvist, in  \cite{FGH}.  This notion of type coincides with the Cohen-Macaulay type of the numerical semigroup ring $\KK[S]$, see \cite{Stamate-2018} for a detailed proof. 

 In general case that $d\geq 1$, the pseudo-Frobenius elements of $S$ are defined by analogy, in  \cite{GORV}, to be $\PF(S)=\{\a\in\NN^d\setminus S \ ; \ \a+S\setminus\{0\}\subseteq S\}$.  As $\NN^d\setminus S$ is not necessarily a finite set, $\PF(S)$ might be an empty set. Indeed $\PF(S)\neq\emptyset$ if and only if $\depth\KK[S]=1$, \cite[Theorem~6]{GORV}.	
 The pseudo-Frobenius numbers of a numerical semigroup,  can be described in terms of Ap\'ery sets. Let $\a$ be a nonzero element of a numerical semigroup $S$. Then $\PF(S)=\{\w-\a \ ; \ \w\in\Max_{\preceq_S}\Ap(S,\a)\}$, see \cite[2.20]{RG}. First, we present a generalization of $\PF(S)$, for simplicial affine semigroups $S\subseteq\NN^d$, in terms of some Ap\'ery sets. 	
 Let $E=\{\a_1,\dots,\a_d\}$ and let $l_i$ be the smallest positive integer such that $l_i\a_{d+i}\in \sum^d_{j=1}\NN\a_j$, for $i=1,\dots,r$. Then 
\[
\Ap(S,E)=\bigcap^d_{i=1}\Ap(S,\a_i)\subseteq\{\sum^r_{i=1}n_i\a_{d+i} \; ; \; 0\leq n_i<l_i\},
\]
is a finite set. The last set in the above equation is  called  $\Gamma$ in \cite{RG-1998}.

\begin{definition}
	The element $\b-\sum^d_{i=1}\a_i$, where $\b\in\Max_{\preceq_S}\Ap(S,E)$,  is called a {\em quasi-Frobenius} element. The set of quasi-Frobenius elements of $S$ is denoted by $\QF(S)$. The number of quasi-Frobenius elements, is called the {\em type} of $S$ and is denoted by $\typ(S)$.
\end{definition}

\begin{remark}\label{rem:QFPF}
 Let $d>1$. If $\f\in \QF(S)\cap\PF(S)$, then $\f+\a_1=\m-\sum_{i=2}^{d}\a_i$, where $\m\in \Max_{\preceq_S}\AP(S,E)$. Since $\f\in \PF(S)$, this follows $\f+\a_1\in S$, which contradicts $\m-\sum_{i=2}^{d}\a_i\notin S$.
\end{remark}

The type of a $d$-dimensional Cohen-Macaulay local ring $(R,\fm)$ is $\typ(R)=\dim_{R/\fm}\Ext^d_R(R/\fm,R)$. 
For a  Cohen-Macaulay ring $R$, the type is defined  as 
the maximum of $\typ(R_\fp)$, where
	$\fp$ ranges in the set of maximal ideals of $R$.

The ring $\KK[S]$ is $\NN$-graded by setting $\deg(\x^\a)=|\a|$, for all $\a\in S$, where $|(a_1,\dots,a_d)|=\sum^d_{i=1}a_i$, denotes the total degree. Therefore, $$\typ(\KK[S])=\typ(\KK[S]_\fm),$$  by \cite[Theorem]{Aoyama-Goto-1975}.  The following result can be also deduced from \cite[Theorem~4.2(ii)]{Campillo-Gimenez-2000}, where the authors provide a combinatorial method to study some simplicial complexes associated to the elements of $S$. We bring a different algebraic proof here.

\begin{proposition}\label{type}
	If $\KK[S]$ is a Cohen-Macaulay ring, then $\typ(S)=\typ(\KK[S]_\fm)=\typ(\KK[\![S]\!])$.
\end{proposition}
\begin{proof}
The ring map $\KK[S]_\fm \longrightarrow \KK[\![S]\!]$ is flat and has only one trivial fiber which is the field $\KK$.  Since $\KK[\![S]\!]/\fm\KK[\![S]\!]$, as an $\KK[\![S]\!]$-module,  has type 1 and depth 0, we have 
\[
\typ(\KK[\![S]\!]\underset{\KK[S]_\fm}{\otimes}\KK[S]_\fm)=\typ(\KK[S]_\fm),
\]	
\[
\depth(\KK[\![S]\!]\underset{\KK[S]_\fm}{\otimes}\KK[S]_\fm)=\depth(\KK[S]_\fm),
\] 
by \cite[Proposition 1.2.16]{Bruns-Herzog}. 
 Thus, $\KK[\![S]\!]$ is Cohen-Macaulay and  $\typ(\KK[S]_\fm)=\typ(\KK[\![S]\!])$.  Let $R=\KK[\![S]\!]$. Then $R$ is a local ring with maximal ideal  $\fm=(\x^{\a_1},\dots,\x^{\a_{d+r}})$. Note that  $\fq=(\x^{\a_1},\ldots,\x^{\a_d})$ is a parameter ideal of $R$, 	since $S$ is simplicial. As $R$ is Cohen-Macaulay,  $\x^{\a_1},\ldots,\x^{\a_d}$ provide  a maximal $R$-regular sequence. By \cite[Lemma~1.2.19]{Bruns-Herzog}, 
\[
\typ(R)=\dim_{R/\fm}(\Hom_R(R/\fm,R/\fq)).
\]
Since $\Hom_R(R/\fm,R/\fq)\cong(0:_{R/\fq}\fm)=\{r\in R/\fq \ ; \ r\fm=0\}$, it is enough to  show that $(0:_{R/\fq }\fm)$
is the $R/\fm$-vector space   generated by residue classes of $\x^{\bs}$, where $\bs\in\max_{\preceq_S}\AP(S,E)$.
For an element $\f\in R$, the residue of $\f$ in $R/\fq$ is equal to the residue of $\sum_{i\geq 1} r_i\x^{\bs_i}$, for some $r_i\in\KK$ and $\bs_i\in \AP(S,E)$.  If the residue of $\f$ in $R/\fq$, belongs to $(0:_{R/\fq}\fm)$, then we derive $\x^{\bs_i+\a_j}\in \fq $, for  $i\geq 1$ and $1\leq j\leq d+r$
which implies $\bs_i\in\max_{\preceq_S}\AP(S,E)$.
Conversely, let $\bs\in\max_{\preceq_S}\AP(S,E)$. Since $\s+\a_i\notin\Ap(S,E)$, for $i=d+1,\dots,d+r$, we get  $\x^{\s+\a_i}\in\fq R$. 
\end{proof}
   	Recall that a Cohen-Macaulay ring is Gorenstein precisely when its Cohen-Macaulay  type is  one. As an immediate consequence of Proposition~\ref{type}, we derive 
   	the following
      \begin{corollary}\cite[4.6 and 4.8]{RG-1998} \cite[4.2]{Campillo-Gimenez-2000}
 $\KK[S]$ is a  Gorenstein ring if and only  if it is Cohen-Macaulay and $\Ap(S,E)$ has a single maximal element with respect to $\preceq_S$.
\end{corollary} 
If $S$ is a numerical semigroup of embedding dimension three, then $\typ(S)\leq 2$, \cite[Theorem~11]{FGH}. The following result is a generalization of  this fact, to simplicial affine semigroups of embedding dimension $d+2$.

\begin{theorem}\label{aa}
	Let  $S$ be of embedding dimension $d+2$. If $\KK[S]$ is Cohen-Macaulay, then $\typ(S)\leq 2$.
\end{theorem}
\begin{proof}
	Let $\Max_{\preceq_S}\Ap(S,E)=\{\m_1,\dots,\m_t\}$. By Lemma~\ref{edim2},  $\m_i$ has a unique maximal expression  $$\m_i=r_{i_1}\a_{d+1}+r_{i_2}\a_{d+2},$$
	for   $i=1,\dots,t$.  Let  $s,k$ be such that $r_{s_1}=\max\{r_{i_1} \ ; \ i=1,\dots,t\}$ and $r_{k_2}=\max\{r_{i_2} \ ; \ i=1,\dots,t\}$.
	Assume on the contrary, that $t\geq 3$ and let $l\in\{1,\dots,t\}\setminus\{s,k\}$. Since $\m_l\npreceq_ S\m_s$ and $\m_l\npreceq_ S\m_k$, we have 
	\[
	r_{k_1}<r_{l_1}<r_{s_1}, r_{s_2}<r_{l_2}<r_{k_2}.
	\]
Moreover,
	\begin{eqnarray}\label{a}
	\m_l+\a_{d+1}=(1+r_{l_1})\a_{d+1}+r_{l_2}\a_{d+2}=\w_1+\sum^d_{i=1}n_i\a_i,\\\m_l+\a_{d+2}=r_{l_1}\a_{d+1}+(1+r_{l_2})\a_{d+2}=\w_2+\sum^d_{i=1}m_i\a_i,\label{a1}
	\end{eqnarray}
	where $\w_1,\w_2\in\Ap(S,E)$ and $n_i,m_i\in\NN$  for $i=1,\dots,d$. Note that 
	\[
\w_1=c_1\a_{d+1}+c_2\a_{d+2} \ , \  \w_2=e_1\a_{d+1}+e_2\a_{d+2},
	\]
for some 	  $c_1,c_2,e_1,e_2\in\NN$. Since $\m_l+\a_{d+i}\notin\Ap(S,E)$, for $i=1,2$, both  $\sum^d_{i=1}n_i\a_i$ and $\sum^d_{i=1}m_i\a_i$ are nonzero. 
As  $r_{l_1}\a_{d+1}+r_{l_2}\a_{d+2}\in \Ap(S,E)$, we get 
	$c_1=e_2=0$. 
	Note that, $1+r_{l_1}\leq r_{s_1}$ and $1+r_{l_2}\leq r_{k_2}$. Therefore, $(1+r_{l_1})\a_{d+1}$ and $(1+r_{l_2})\a_{d+2}$  belong to $\AP(S,E)$. Consequently, 
	\begin{equation}\label{res}
	r_{l_1}>e_1 \text{ and }  r_{l_2}>c_2,
	\end{equation}
	which implies 
	$(1+e_1)\a_{d+1}$ and $(1+c_2)\a_{d+2}$ are also in  $\Ap(S,E)$.  
	
	Now,  subtracting  (\ref{a1}) from (\ref{a}), we derive
	\[
	(1+e_1)\a_{d+1}-(1+c_2)\a_{d+2}=\sum^d_{i=1}(n_i-m_i)\a_i.
	\]
	Since  $\KK[S]$ is Cohen-Macaulay,  Proposition~\ref{CM} implies that 
	\[
	(1+e_1)\a_{d+1}=(1+c_2)\a_{d+2}.
	\]
	Without loss of generality we may assume that $e_1>c_2$. 	Then
	\[
	r_{l_2}\a_{d+2}=(r_{l_2}-c_2-1)\a_{d+2}+(1+e_1)\a_{d+1}.
	\]  Consequently, $\ord(r_{l_2}\a_{d+2})\geq r_{l_2}-c_2+e_1>r_{l_2}$, in contradiction to our choice of maximal expression of $\m_l$.
\end{proof}

The converse of Theorem~\ref{aa} is not true.
\begin{example}
	Let $\a_1=(2,0),\a_2=(0,2),\a_3=(4,1),\a_4=(2,3)$. Since $\a_1-\a_2=\a_3-\a_4$, $\KK[S]$ is not Cohen-Macaulay, by Proposition~\ref{CM}. However,
	$$\Ap(S,E)\setminus\{0\}=\max_{\preceq_S}\Ap(S,E)=\{(4,1), (2,3)\}.$$
	
\end{example}

	The following two  examples show that either if  $\KK[S]$ is not Cohen-Macaulay or if $r\geq 3$, then $\typ(S)$ does not have any upper bound in terms of its embedding dimension.

\begin{example}\label{3^n}
Let $\a_1=(3,0),\a_2=(0,3^n),\a_3=(5,2),\a_4=(2,2+3^n)$, where $n$ is a positive integer. Since $\a_4-\a_3=\a_2-\a_1$, $\KK[S]$ is not Cohen-Macaulay, by Proposition~\ref{CM}. 
First, we show that 
\begin{equation}\label{3.4}
\Ap(S,E)=\{r\a_3+s\a_4 ;  \text{ for all }  r,s\in \NN \ \text{such that} \ r+s<3^n \}.
\end{equation}
Assume on the contrary that, $r\a_3+s\a_4\notin\AP(S,E)$
for some $r,s\in \NN$ such that $r+s<3^n$. Then $r\a_3+s\a_4=\sum^4_{i=1}n_i\a_i$ for some integers  $n_i\geq 0$ such that either $n_1>0$ or $n_2>0$. Consequently
\[
(r-n_3)(5,2)+(s-n_4)(2,2+3^n)=n_1(3,0)+n_2(0,3^n),
\]
which implies 
\begin{eqnarray}
2(r+s-n_3-n_4)+3^n(s-n_4)=3^nn_2 \label{3.5}\\
5(r-n_3)+2(s-n_4)=3n_1. \label{3.6}
\end{eqnarray}
Therefore, $r+s=3^nk+n_3+n_4$ for some integer $k$. Since $r+s<3^n$, we get $k\leq 0$. Then $r+s-n_3-n_4\leq 0$. Since $n_2\geq 0$, we get by (\ref{3.5}) that   $s-n_4\geq0$ and $r-n_3\leq0$. Note that
\[
3(r-n_3)+2(r+s-n_3-n_4)=3n_1,
\]
from (\ref{3.6}). Since $r-n_3\leq0$, $r+s-n_3-n_4\leq 0$ and $n_1\geq 0$, we get 
\[
r-n_3=r+s-n_3-n_4=n_1=0.
\]
Then $n_2=0$, from (\ref{3.5}), a contradiction.
We have shown that,  the set on the right hand side of (\ref{3.4}) is a subset of $\Ap(S,E)$.
Now, let $r,s\in\NN$ such that  $r+s=3^n.$ So
\begin{align*}
r\a_3+s\a_4-\a_1&=r(5,2)+s(2,2+3^n)-(3,0)
\\&
=(2(r+s)+3(r-1),2(r+s)+3^{n}s)
\\&
=( 2\times3^{n-1}+(r-1))(3,0)+(2+s)(0,3^n).
\\&
\end{align*}
Therefore $r\a_3+s\a_4\notin \AP(S,E)$. If $r+s>3^n$, then
\[
 r\a_3+s\a_4=r^{\prime}\a_3+s^{\prime}\a_4+(r-r^{\prime})\a_3+(s-s^{\prime})\a_4
 \]
 for some nonnegative integers $r',s'$ such that $r>r'$, $s>s'$ and  $r^{\prime}+s^{\prime}=3^n$.   
Consequently, $r\a_3+s\a_4\notin \AP(S,E)$. 
Therefore, (\ref{3.4}) is proved.  Now, one can easily see that 
$$\max_{\preceq_S}\Ap(S,E)=\{r\a_3+s\a_4 ;  \text{ for all }  r,s\in \NN \ \text{such that} \ r+s=3^n-1\},$$
and so $\typ(S)=|\{(r,s)\in\NN^2 \ ; \ r+s=3^n-1\}|=3^n$.
\end{example}

\begin{example}\label{ex:3.7}
	For an integer $a\geq 3$, let $S$ be the affine semigroup generated by  $\a_1=(a^2,0), \a_2=(0,a^2),\a_3=(a^2-a,a^2-a),\a_4=(a^2-a+1,a^2-a+1),\a_5=(a^2-1,a^2-1)$. Then $S$ is simplicial with extremal rays $\a_1,\a_2$. Let $T$ be the numerical semigroup generated by $\{a^2-a, a^2-a+1, a^2-1, a^2 \}$. Then
\begin{equation}\label{AP}
	\AP(S,E)=\AP(S,\a_1+\a_2)=\{(s,s) \ ; s\in\AP(T,a^2)\}.
\end{equation}
	Therefore, $\typ(S)=\typ(T)=2a-4$, by \cite[(3.4)Proposition]{Cavaliere-Niesi-1983}. 
Note that $\KK[S]$ is  Cohen-Macaulay, by Lemma~\ref{rmk:cm}.
\end{example}

\section{The conductor ideal of simplicial affine semigroups}\label{sec:norm}
 Thorough this section, 
$S\subseteq\NN^d$  is a  simplicial affine semigroup with $\mgs(S)=\{\a_1,\dots,\a_{d+r}\}$,  where  $\a_1,\dots,\a_d$ are the  extremal rays of $S$. Let $R=\KK[S]$ be the affine semigroup ring. 
 Recall that the normalization of an integral domain $R$ is the set of elements
	in its field  of fractions satisfying a monic polynomial in $R[y]$. Then 
 $R=\KK[S]$ is an integral domain with  normalization $\bar{R}=\KK[\Gr(S)\cap\co(S)]$
 \cite[Proposition~7.25]{Miller-Sturmfels}.
 Recall that the  conductor of $R$,  $C_R=(R:_T\bar{R})$,  where $T$ denotes the total ring of fractions of $R$, is the largest common ideal of $R$ and $\bar{R}$, \cite[Exercise~2.11]{Huneke-Swanson}. 
The purpose of this section, is to investigate the normality  of  $R$ and the generating set of $C_R$ as an ideal of $\bar{R}$. 

 The integral closure of $S$ in $\Gr(S)$, $\bar{S}=\{\a\in\gr(S) \ ;\ n\a\in S \text{ for some } n\in\NN\}$, is called the {\em normalization} of $S$. 
   As a geometrical interpretation, one can see that  $\bar{S}=\co(S)\cap\gr(S)$.  The semigroup  $S$ is {\em normal} when $S=\bar{S}$, equivalently $\KK[S]$ is a normal ring, \cite{Bruns-Gubeladze-2009,Bruns-Herzog}.   Since $S$ is finitely generated, $\co(S)$ is generated by  finitely many rational vectors, i.e. it is the intersection of
 finitely many rational vector halfspaces, \cite[Corollary~7.1(a)]{Schrijver}. By Gordan's lemma, $\bar{S}$ is also finitely generated. 

The {\em conductor} of $S$ is defined as $\fc(S)=\{\b\in S \ ; \ \b+\bar{S}\subseteq S\}$. The conductor, $\fc(S)$, is   the largest ideal of $S$ that is also an ideal of $\OS$, \cite[Exercise~2.9]{Bruns-Gubeladze-2009}. 

	\begin{remark}\label{rem:normal}
		As $\fc(S)$ is an ideal of $S$, we have $S=\bar{S}$ precisely when $0\in\fc(S)$. In other words, $S$ is normal if and only if $\fc(S)=S$.
	\end{remark}

When $S$ is fully embedded in $\NN^d$,  that is the affine
 subspace it generates coincides with $\RR^d$, we have $\Gr(S) \cong\ZZ^d$. In the case that   $\Gr(S)=\ZZ^d$, we have $\bar{S}=\NN^d\cap\co(S)$.  
  The later property may happen also for affine semigroups that $\Gr(S)\neq\ZZ^d$.  For instance, it holds when $({\co(S)}\setminus S)\bigcap \NN^d$ is finite, such semigroups are considered in \cite{C-sem}.
 	
 	\begin{lemma}
 		If  $({\co(S)}\setminus S)\bigcap \NN^d$  is finite, then $\gr(S)\cap\co(S)=\NN^d\cap\co(S)$. In particular, $S$ is normal if and only if $S=\NN^d\cap\co(S)$. 
 	\end{lemma}
 	\begin{proof}
 All vectors in $\co(S)$ have nonnegative components. Thus $\gr(S)\cap\co(S)\subseteq\NN^d\cap\co(S)$. 
 		Let $\a\in \NN^d\cap\co(S)$. Then $\a+\sum^d_{i=1}l_i\a_i\in\co(S)\cap\NN^d$, for all $l_1,\dots,l_d\in\NN$. Since $(\co(S)\cap\NN^d)\setminus S$ is a finite set, we have 
 		$\a+\sum^d_{i=1}l_i\a_i\in S$, for some $l_1,\dots,l_d\in\NN$, which implies $\a\in\Gr(S)$.
 	Therefore, $\NN^d\cap\co(S)=\gr(S)\cap\co(S)$. 
 	\end{proof}

 	\begin{lemma}\label{lem:OS}
 			As an affine semigroup, $\bar{S}$ is generated by $(P_S\cap\Gr(S))\cup\{\a_1,\dots,\a_d\}$,  and
$P_S\cap \Gr(S)=\{r(\w) \ ;\  \w\in \AP(S,E)\}$.
	\end{lemma} 

\begin{proof}
An element $\a\in\bar{S}$ can be written as $\a=r(\a)+\sum^d_{i=1}r_i\a_i$	for some $r_1,\dots,r_d\in~\NN$.
Therefore, $(P_S\cap\Gr(S))\cup\{\a_1,\dots,\a_d\}$ provides a generating set for $\bar{S}$.
 For the last statement, let $\b\in P_S\cap\Gr(S)$. By Remark~\ref{rmk:GrAp},  there exist $r_1,\dots,r_d,s_1,\dots,s_d\in\NN$ and $\w\in\Ap(S,E)$ such that   $\b+\sum^d_{i=1}r_i\a_i=\w+\sum^d_{i=1}s_i\a_i$.  Since $\b\in P_S$, it follows that $\b=r(\w)$. 
\end{proof}	 
As an immediate consequence of   Lemma~\ref{lem:OS}, 
 	\begin{equation}\label{eq:c(S)}
 	\fc(S)=\{\a\in S \ ; \ \a+r(\w)\in S \text{ for all }  \w\in\Ap(S,E) \}.
 	\end{equation}

\begin{example}\label{F+1}
	Let $S$ be a numerical semigroup, that is  $d=1$. Then $\OS=\NN$ and $\fc(S)=\{n \ ; \ n\geq F+1\}$, where $F$ is the maximal integer in $\ZZ\setminus S$.  Therefore, $\fc(S)$  is  generated by $\{F+1\}$ as an ideal of $\NN$. In this regard, $F+1$ is called the conductor of $S$. 
\end{example}

\begin{lemma}\label{4.5}
		The following statements are equivalent.
	\begin{enumerate}
		\item $-\QF(S)\subseteq\co(S)$;
		\item  $\Ap(S,E)\subseteq \bar{P}_S$.
	\end{enumerate}
\end{lemma}
\begin{proof}
	(1)$\implies$(2): 
Let $\w\in \AP(S,E).$ Then $\m-\w\in S$, for some $\m \in \Max_{\preceq_S}\AP(S,E)$.   As $\m-\sum^d_{i=1}\a_i\in\QF(S)$, we get  $\sum_{i=1}^{d}(1-[\m]_i)\a_i\in \co(S)$ precisely when $[\w]_i\leq[\m]_i\leq1$.

(2)$\implies$(1):  Let $\f=\m-\sum^d_{i=1}\a_i\in\QF(S)$ for some $\m \in \Max_{\preceq_S}\AP(S,E)$. As  $[\m]_i\leq 1$ for $i=1,\dots,d$, we get $-\f=\sum_{i=1}^{d}(1-[\m]_i)\a_i\in\co(S)$.	
\end{proof}	

If we replace $\bar{P}_S$ in the above lemma, with $P_S$, then we derive an equivalent condition for $S$ to be normal. 
Roughly speaking,  having more negative coefficients $[\f]_i$ for $\f\in \QF(S)$, makes the semigroup more close to being normal. 
In this order,  we need to consider the relative interior of $\co(S)$. 
  Let $\relint(S)$ denote the elements of $\RR^d$ that belong to the relative interior of $\co(S)$,  
	\[
	\relint(S)=\left\{\b\in \co(S) \ ; \ \b=\sum_{i=1}^d \lambda_i \a_i \text{ with  $\lambda_i \in \RR_{>0}$ for all $i=1, \dots, d$}\right\}.
	\]
 
\begin{theorem}\label{QF&NORMAL}
	The following statements are equivalent.
	\begin{enumerate}
		\item $S$ is normal;
		\item $-\QF(S)\subseteq S\cap\relint(S)$;
		\item $-\QF(S)\subseteq\relint(S)$;
		\item  $\Ap(S,E)\subseteq P_S$.
		\end{enumerate}
\end{theorem}

\begin{proof} (1)$\implies$(2): Let $\f\in\QF(S)$. Then $\f=\m-\sum^d_{i=1}\a_i$, for some $\m\in\Max_{\preceq_S}\Ap(S,E)$.
	If $[\m]_j\geq 1$, for some $1\leq j\leq d$, then     $\m=\a_{j}+([\m]_j-1)\a_j+\sum^d_{i=1,i\neq j}[\m]_i\a_i$, which implies   $\m-\a_j\in\bar{S}=S$, which contradicts   $\m\in\Ap(S,E)$. Therefore $[\m]_i<1$, for $i=1,\dots,d$. 
	Now, $-\f=\sum^d_{i=1}(1-[\m]_i)\a_i\in\relint(S)\cap\Gr(S)\subseteq S$. 

(2)$\implies$(3) is clear.

	(3)$\implies$(4): Let $\w\in \AP(S,E).$ Then $\m-\w\in S$, for some $\m \in \Max_{\preceq_S}\AP(S,E)$.   As $\m-\sum^d_{i=1}\a_i\in\QF(S)$, we get  $\sum_{i=1}^{d}(1-[\m]_i)\a_i\in \relint(S)$.  Thus, $1-[\m]_i>0$ which implies  $[\w]_i\leq[\m]_i<1$. 

(4)$\implies$(1):  $\Ap(S,E)\subseteq P_S$ is equivalent to  $\Ap(S,E)=r(\Ap(S,E))$. As the later equality holds precisely when $r(\Ap(S,E))\subseteq S$, the  result follows by Lemma~\ref{lem:OS}.
\end{proof}

 Our next aim in this section is to find a generating  set for $\fc(S)$ as an ideal of $\OS$.   
  Recall from Section~\ref{sec:Back}, that $C_j=\{\w\in\Ap(S,E) \ ; \ r(\w)=\b_j\}$, for $ j=0,\dots,k$, where 
  $r(\Ap(S,E))=\{r(\w) \ ; \ \w\in\Ap(S,E)\}=\{0=\b_0,\b_1,\dots,\b_k\}.$ For any $(\w_1,\dots,\w_k)\in C_1\times\dots\times C_k$, we consider  the  vector
 	\[
  	\f_{(\w_1,\dots,\w_k)}=\sum^d_{i=1}f_i\a_i,
  	\]
where  $f_i=\max\{[\w_j-r(\w_j)]_i \ ;\  j=1,\dots,k\}$,  for $i=1,\dots,d$.  Note that 
  	$$f_i=\max\{\lfloor[\w_j]_i\rfloor \ ; \ j=1,\dots,k\},$$  for $i=1,\dots,d$,  where $\lfloor[\w_j]_i\rfloor$ denotes the greatest integer less than or equal to $[\w_j]_i$.

   \begin{lemma}\label{lem:A}
  		Given $(\w_1,\dots,\w_k)\in C_1\times\dots\times C_k$, the vector 	$\f_{(\w_1,\dots,\w_k)}$ belongs to $\fc(S)$.
  \end{lemma}
\begin{proof}
Let $\f=\f_{(\w_1,\dots,\w_k)}$. 
By (\ref{eq:c(S)}), it is enough  to show that $\f+\b_j\in S$ for  $j=1,\dots,k$. Let $\b_j=r(\w_j)$. Note that 
$\f=\sum^d_{i=1}([\w_j-\b_j]_i+r_i)\a_i$, for some $r_i\in\NN$. Therefore, $\f+\b_j=\w_j+\sum^d_{i=1}r_i\a_i\in S$, for $j=1,\dots,k$.
\end{proof}

  	\begin{corollary}
 The following statements are equivalent.
 \begin{enumerate}
 	\item $S$ is normal;
 	\item $\f_{(\w_1,\dots,\w_k)}=0$, for all $\w_i\in C_i$ and $i=1,\dots,k$;
 	\item $\f_{(\w_1,\dots,\w_k)}=0$, for some $\w_i\in C_i$ and $i=1,\dots,k$.
 \end{enumerate}	
  	\end{corollary}
  	\begin{proof}
(1)$\implies$(2) follows from Theorem~\ref{QF&NORMAL}. 

(2)$\implies$(3) is clear.

(3)$\implies$(1): By Lemma~\ref{lem:A}, $0=\f_{(\w_1,\dots,\w_k)}$ belongs to $\fc(S)$ which is an ideal of $S$. Therefore, $S$ is normal, see Remark~\ref{rem:normal}.
\end{proof}
  
 \begin{theorem}\label{fff}
Let $\c$ be a minimal generator of $\fc(S)$. Then there exist $(\w_{1},\dots,\w_{k})\in C_1\times\dots\times C_k$ such that $\c=\f_{(\w_{1},\dots,\w_{k})}-\b_j+\sum^d_{i=1}l_i\a_i$ for some $l_i\in\{0,1\}$ and $j\in\{0,\ldots,k\}$.  Moreover, at least for one $i$, we have $l_i=0$.

\end{theorem}
\begin{proof}
 Since $\c+\b_j\in S$, for $j=0,\dots,k$,
\begin{equation}\label{eq:c+bi}
\c+\b_j=\w_{t_j}+\sum^d_{i=1}r_{j_i}\a_i,
\end{equation}	
 for some $\w_{t_j}\in\Ap(S,E)$ and $r_{j_i}\in\NN$.
Note that $r(\w_{t_j})\neq r(\w_{t_i})$, for $0\leq i\neq j\leq k$, since otherwise $\b_j-\b_i\in\Gr(\a_1,\dots,\a_d)$ which is not possible.
	Let $\{t_1,\dots,t_k\}=\{1,\dots,k\}$ such that $\w_{i}\in C_i$, and let  $\f=\f_{(\w_{1},\dots,\w_{k})}$. 
	Then $$[\f]_j=\max\{\lfloor[\w_i]_j\rfloor \ ; \ i=1,\dots,k\},$$ for $j=1,\dots,d$. 
Let $1\leq s\leq d$.
 If $r_{j_s}\geq 1$ for  $j=0,\dots,k$, then 
\[
\c-\a_s+\b_i=\w_{t_i}+\sum^d_{j=1,j\neq s}r_{i_j}\a_j+(r_{i_s}-1)\a_s\in S,
\]
for $i=0,\dots,k$, which implies $\c-\a_s\in\fc(S)$, a contradiction.	Thus, $r_{j_s}=0$ for some $0\leq j\leq k$. Consider  $0\leq h\leq k$ such that 
\[
[\b_h]_s=\max\{[\b_j]_s \ ; \  r_{j_s}=0 \ , \ 0\leq j\leq k\}.
\]
Then 
$[\w_{t_h}]_s=[\c+\b_h]_s=\max\{[\w_{t_j}]_s \ ; \ r_{j_s}=0 \ , \ 0\leq j\leq k\}$.
 If $r_{i_s}>0$ for some $0\leq i\leq k$, then 
\[
[\w_{t_h}]_s-[\b_h]_s=[\c]_s=[\w_{t_i}]_s+r_{i_s}-[\b_i]_s>[\w_{t_i}]_s.
\]
Consequently, 
\[
[\w_{t_h}]_s=\max\{[\w_i]_s \ ; \ i=1,\dots,k\},
\]
and thus 
 \[
\lfloor[\w_{t_h}]_s\rfloor=[\f]_s.
\]
 Let $c_s=\lceil[\c]_s\rceil$, where $\lceil[\c]_s\rceil=\operatorname{ceil}([\c]_s)$ denotes the least integer greater than or equal to $[\c]_s$, for $s=1,\dots,d$.
 Now, we distinguish the following two possibilities:
\begin{enumerate}
	\item If $[\b_h]_s\geq  c_s-[\c]_s$, then  $[\w_{t_h}]_s=[\b_h+\c]_s\geq c_s$. As $c_s$ is an integer, it follows that  $\lfloor[\w_{t_h}]_s\rfloor\geq c_s$. If $\lfloor[\w_{t_h}]_s\rfloor>c_s$, then $\lfloor[\w_{t_h}]_s\rfloor\geq 1+c_s=1+\lceil[\c]_s\rceil$, which implies   $[\b_h]_s=[\w_{t_h}]_s-[\c]_s\geq1$, a contradiction. 
Therefore, $$c_s=\lfloor[\c+\b_h]_s\rfloor=\lfloor[\w_{t_h}]_s\rfloor=[\f]_s.$$
	\item If $[\b_h]_s<c_s-[\c]_s$, then  $\lfloor[\b_h+\c]_s\rfloor=\lfloor[\c]_s\rfloor$.  Note that, $c_s-[\c]_s>0$ which means $[\c]_s$ is not an integer. Thus, $\lfloor[\c]_s\rfloor=c_s-1$. Therefore, 
	\[
	c_s=\lfloor[\c]_s\rfloor+1=\lfloor[\b_h+\c]_s\rfloor+1=\lfloor[\w_{t_h}]_s\rfloor+1=[\f]_s+1.
	\]
\end{enumerate}
If $c_s$ satisfies (1), then let  $l_s=0$, and otherwise let $l_s=1$. Then $$\c=\f-\sum^d_{s=1}(c_s-[\c]_s)\a_s+\sum^d_{i=1}l_i\a_i.$$
Note that $\sum^d_{s=1}(c_s-[\c]_s)\a_s=\sum^d_{s=1}c_s\a_s-\c$ belongs to $\gr(S)\cap P_S=r(\Ap(S,E))$, see Lemma~\ref{lem:OS}. 

 For the last statement, let $1\leq j\leq k$ and let  $$\{i_1,\dots,i_t\}=\{i \ ; \ 1\leq i\leq d,  [\b_j]_i>0 \}.$$ Then $\sum^t_{s=1}\a_{i_s}-\b_j\in \gr(S)\cap P_S=r(\Ap(S,E))$. 
In particular, $\f+\sum^d_{i=1}\a_i-\b_j\in \{\f\}+\bar{S}$.
As $\f\in\fc(S)$ by Lemma~\ref{lem:A}, it means that $\f+\sum^d_{i=1}\a_i-\b_j$ is not in the minimal generating set of $\fc(S)$.
\end{proof}

\begin{corollary}\label{c&a} 
	If $\c\in\fc(S) $ such that $[\c]_i\in \NN$ for $i=1,\dots,d$, then  $[\c]_i\geq[\f_{(\w_1,\dots,\w_k)}]_i$
	for some  $\w_1,\dots,\w_k\in C_1\times\dots\times C_k$.
\end{corollary}
\begin{proof}
	Note that $\c=\c'+\b$ for a minimal generator $\c'$ of $\fc(S)$ and some $\b\in\bar{S}$. 
	By Theorem~\ref{fff}, there exist  $(\w_{1},\dots,\w_{k})\in C_1\times\dots\times C_k$ and $l_1,\dots,l_d\in\{0,1\}$ such that $\c'=\f_{(\w_{1},\dots,\w_{k})}-{\b_j}+\sum^d_{i=1}l_i\a_i$ for some $j\in\{0,\ldots,k\}$. Thus
	\[
	[\c]_i\geq[\f_{(\w_{1},\dots,\w_{k})}]_i-[\b_j]_i.
	\]
	As $[\c]_i$ is an integer, and $0\leq [\b_j]_i<1$, we have $[\c]_i\geq [\f_{(\w_{1},\dots,\w_{k})}]_i$.
\end{proof}

 \begin{example}\label{ex:4.7}
	Let $\a_1=(3,0), \a_2=(0,3), \a_3=(5,2), \a_4=(2,5)$. As we have seen in Example~\ref{3^n},
	\begin{eqnarray*}
		\AP(S,E)&=&\{0,\a_3,\a_4,\a_3+\a_4,2\a_3,2\a_4\}\\
		&=&\{0,\w_1=(5,2),\w_2=(2,5),\w_3=(7,7),\w_4=(10,4),\w_5=(4,10)\}
	\end{eqnarray*}
	and $r(\Ap(S,E))=\{0,\b_1=(1,1),\b_2=(2,2)\}$. Note that $C_1=\{\w_3,\w_4,\w_5\}$, $C_2=\{\w_1,\w_2\}$  and 
	$\f_{(\w_3,\w_i)}=2\a_1+2\a_2=(6,6)$, $\f_{(\w_4,\w_i)}=3\a_1+\a_2=(9,3)$, $\f_{(\w_5,\w_i)}=\a_1+3\a_2=(3,9)$, for $i=1,2$. 

As $\{(6,6)-(1,1),  (9,3)-(1,1), (3,9)-(1,1)\}+r(\Ap(S,E))\subset S$, we have $$\{(5,5),(8,2),(2,8)\}\subset\fc(S).$$
If $\fc(S)\neq \{(5,5),(8,2),(2,8)\}+\bar{S}$,  the other generators of $\fc(S)$ are among 
\[
\{(9,3), (3,9), (6,6)\}+\{l_i\a_i-(2,2) \ ; \ l_i\in\{0,1\} , i=1,2 \},
\]
by Theorem~\ref{fff}.
Since the above set which equals $$\{(7,1), (1,7), (4,4), (10,1), (1,10), (4,7), (7,4)\},$$
has no element in $S$,   $\fc(S)$ is generated by $\{(5,5), (8,2 ), (2,8)\}=\{\f_{(\w_3,\w_1)}-\b_1, \f_{(\w_4,\w_1)}-\b_1, \f_{(\w_5,\w_1)}-\b_1\}$, as an ideal of $\bar{S}=\langle(3,0),(0,3),(1,1)\rangle$.
\end{example}

 The following example shows that the summand $\sum^d_{i=1}l_i\a_i$ in the statement of Theorem~\ref{fff} can not be removed. 

\begin{example}\label{New-Ex}
	Let $\a_1=(5,2), \a_2=(2,2), \a_3=(2,1),\a_4=(5,3)$. Then $\Ap(S,E)=\{0, \w_1=(2,1), \w_2=(4,2), \w_3=(6,3), \w_4=(8,4), \w_5=(5,3)\}$ and $r(\Ap(S,E))=\{0, \b_1=(2,1), \b_2=(4,2), \b_3=(1,1), \b_4=(3,2), \b_5=(5,3)\}$. Note that $C_i=\{\w_i\}$ for $i=1,\dots, 5$ and $\f_{(\w_1,\dots,\w_5)}=\a_1$. 
	By Theorem~\ref{fff}, the generators of $\fc(S)$ are among
\[
\{\a_1-\b_i , 2\a_1-\b_i, \a_1+\a_2-\b_i \ ; \ i=0,\dots,5\}.
\]	
The only elements of the above set, that belong also to $S$ are 
$$\{(5,2), (10,4), (5,3), (2,1), (7,4), (4,2), (6,3)\}.$$
 Note that $(2,1)+(1,1)\notin S$,  $\{(5,2),  (4,2)\}+r(\Ap(S,E))\subseteq S$, $\{(10,4), (7,4), (6,3)\}\subset (5,2)+\bar{S}$ and $(5,3)=(4,2)+\b_3$. Therefore,  $\fc(S)$ is generated by 
 $\{(5,2), (4,2)\}=\{\f_{(\w_1,\dots,\w_5)},  \f_{(\w_1,\dots,\w_5)}+\a_2-\b_4\}$, as an ideal of $\bar{S}=\langle (1,1),(2,1),(5,2)\rangle$.
\end{example}

\begin{proposition}\label{prop:sing-max}
Assume that there is  a fixed class $C_j$ such that for any $\w\in C_j$ and $\w'\in\Ap(S,E)\setminus C_j$, one has  $\max_{\preceq_c}(\w,\w')=\w$. If either $C_j$  is a singleton or $\b_j=\min_{\preceq_{c}}(r(\Ap(S,E))\setminus\{0\})$, then  $\fc(S)$ is generated by 
	$$\{\w-\b \ ; \ \w\in C_j \ , \ \b\in\max_{\preceq_c}\{\b_1,\dots,\b_k\}\},$$ as an ideal of $\bar{S}$.
\end{proposition}

\begin{proof}
    First we show  that
	\begin{equation}\label{ij-ls}
	\w-\b_l+\b_s\in S,
	\end{equation}
	for any $\w\in C_j$ and  $0\leq l,s\leq k$.
	
As	$\w-\b_l+\b_s\in\Gr(S)$, we have by Remark~\ref{rmk:GrAp}, that $\w-\b_l+\b_s+\sum^d_{i=1}r_i\a_i=\w'+\sum^d_{i=1}s_i\a_i$ for some $\w'\in\Ap(S,E)$ and nonnegative integers  $r_1,\dots,r_d,s_1,\dots,s_d$. If $r_1=\dots=r_d=0$, then (\ref{ij-ls}) is clear. Assume that $r_h>0$, for some $1\leq h\leq d$. Then $s_h=0$ by our choice in Remark~\ref{rmk:GrAp}, and  $[\w]_h<[\w']_h$.  If $\w'\notin C_j$, then $\max_{\preceq_c}(\w,\w')=\w$ which implies $[\w]_h\geq [\w']_h$,  a contradiction.
 Thus  $\w'\in C_j$. In other words, $\w$ and $\w'$ have the same remainder $\b_j$. Therefore, $\b_s=\b_l$, and (\ref{ij-ls}) holds. Consequently,   $\w-\b_r\in\fc(S)$ for $r=0,\dots,k$.

Let $\c$ be a  minimal generator  of $\fc(S)$. By Theorem~\ref{fff},  $\c$  can be written as 
\begin{equation}\label{4.4}
\c=\f_{(\w_{1},\dots,\w_{k})}-\b_t+\sum^d_{i=1}l_i\a_i,
\end{equation}
 for some $(\w_1,\dots,\w_k)\in C_1\times\dots\times C_k$,  $0\leq t\leq k$ and some $l_1,\dots,l_d\in\{0,1\}$. Note that $\f_{(\w_{1},\dots,\w_{k})}=\w_j-\b_j$ and $\b_j+\b_t=\b_s+\sum^d_{i=1}l'_i\a_i$ for some  $l'_1,\dots,l'_d\in\{0,1\}$ and $0\leq s\leq k$. 
 Therefore, we get 
\begin{equation}\label{4.5}
\c=\w_j-(\b_j+\b_t)+\sum^d_{i=1}l_i\a_i=\w_j-\b_s+\sum^d_{i=1}(l_i-l'_i)\a_i,
\end{equation}
from (\ref{4.4}). If $\b_s=0$, then $\c$ and $\w_j$ are in the same congruence class modulo the group spanned by the extremal rays. In other words, $r(\c)=r(\w_j)=\b_j$. Thus, $\c=\w+\sum^d_{i=1}h_i\a_i$, for some $\w\in C_j$ and $h_1,\dots,h_d\in\NN$. This provides a contradiction with minimality of $\c$, as $\w-\b_i\in\fc(S)$, for $1\leq i\leq k$. Therefore, $\b_s\neq0$.  Now, we distinguish the following two cases:

Case 1, $\b_j=\min_{\preceq_c}(r(\Ap(S,E))\setminus\{0\})$. Then $[\b_j]_i\leq[\b_s]_i=[\b_j]_i+[\b_t]_i-l'_i$, which implies $l'_i=0$, for $i=1,\dots,d$.
 
Case 2:  $C_j$  is a singleton. From equation (\ref{4.5}), we have  $r(\c+\b_s)=r(\w_j)=\b_j$ which implies  $\c+\b_s=\w_j+\sum^d_{i=1}r_i\a_i$ for some $r_1,\dots,r_d\in\NN$. Now, looking again at (\ref{4.5}), we derive 
$$\w_j+\sum^d_{i=1}r_i\a_i=\w_j+\sum^d_{i=1}(l_i-l'_i)\a_i.$$
As $\a_1,\dots,\a_d$ are linearly independent,  we get $l_i-l'_i=r_i\geq0$, for $i=1,\dots,d$.

Thus in both cases, $l_i-l'_i\geq0$, for $i=1,\dots,d$.  Since $\w_j-\b_s\in\fc(S)$ and $\c$ is a minimal generator of $\fc(S)$, we derive from (\ref{4.5}), that $l_i-l'_i=0$ for $i=1,\dots,d$ and  $\c=\w_j-\b_s$.

If $\b_s\preceq_c\b_r$ for some $1\leq r\leq k$, then $\b_r-\b_s\in\co(S)\cap\gr(S)=\bar{S}$ and $\w_{j}-\b_s=\w_{j}-\b_r+\b_r-\b_s$. As $\c$ is a minimal generator, we get $\b_s=\b_r$. Thus,  $\b_s\in\max_{\preceq_c}\{\b_1,\dots,\b_k\}$. 
	\end{proof}

 Applying the above proposition to the semigroup in Example~\ref{ex:4.7}, provides an easier argument to find the minimal generating set of $\fc(S)$.   
 \begin{example}
 	Let $\a_1=(3,0), \a_2=(0,3), \a_3=(5,2), \a_4=(2,5)$. As we have seen in  Example~\ref{ex:4.7},  $\AP(S,E)=\{0,\w_1=(5,2),\w_2=(2,5),\w_3=(7,7),\w_4=(10,4),\w_5=(4,10)\}$, $r(\Ap(S,E))=\{0,\b_1=(1,1),\b_2=(2,2)\}$,  $C_1=\{\w_3,\w_4,\w_5\}$ and  $C_2=\{\w_1,\w_2\}$. Note that
 $\max_{\preceq_c}\{\w_i,\w_j\}=\w_j$ for $i=1,2$ and $j=3,4,5$.  	
	Therefore, $\fc(S)$ is generated by  $\{\w_3-(2,2),\w_4-(2,2),\w_5-(2,2)\}=\{(5,5),(8,2),(2,8)\}$, as an ideal of $\bar{S}=\langle(3,0),(0,3),(1,1)\rangle$.
\end{example}

\begin{example}
	 Let $\a_1=(2,0),\a_2=(0,2),\a_3=(4,1),\a_4=(2,3).$ We have $\AP(S,E)=\{0,\w_1=(4,1),\w_2=(2,3)\}$ and $r(\Ap(S,E))=\{0,(0,1)\}$. More precisely,
$\w_1-2\a_1=\w_2-\a_1-\a_2=1/2\a_2$. Thus  $k=1$ and $C_1=\{(4,1),(2,3)\}$. By Proposition~\ref{prop:sing-max}, $\fc(S)$ is generated by $\{(4,1)-(0,1),(2,3)-(0,1)\}=\{(4,0),(2,2)\}$ as an ideal of $\bar{S}=\langle (2,0),(0,1)\rangle$.
\end{example}

\begin{remark}
If $C_i$ is a singleton for $i=1,\dots,k$, then the hypothesis on the existence of $j$ in  Proposition~\ref{prop:sing-max} is equivalent to existence of a single maximal element of $\Ap(S,E)$ with respect to $\preceq_c$. This condition  can not be removed. For instance, let $S$ be the semigroup defined in Example~\ref{New-Ex}. Then $\max_{\preceq_c}\Ap(S,E)=\{\w_4,\w_5\}=\{(8,4), (5,3)\}$ and $\max_{\preceq_c}(r(\Ap(S,E)))=\{\b_5=(5,3)\}$. But, $(8,4)-(5,3)$ is not in $S$ and in particular, it is not in $\fc(S)$.
\end{remark}

 \begin{corollary}\label{cor:principal}
	If  $\KK[S]$ is a Cohen-Macaulay ring,  and $\max_{\preceq_c}\Ap(S,E)=\{\w\}$, then  $\fc(S)$ is generated by 
	$\{\w-\b \ ; \b\in\max_{\preceq_c}r(\Ap(S,E))\},$ as an ideal of $\bar{S}$. In particular, if    $\max_{\preceq_c}(r(\Ap(S,E)))=\{\b\}$, then $\fc(S)$ is a principal ideal of $\bar{S}$ generated by $\w-\b$. 
\end{corollary}

As Example~\ref{New-Ex} shows, for an  affine semigroup with  Cohen-Macaulay semigroup ring,  $\max_{\preceq_c}\Ap(S,E)$ is not necessarily a singleton. Here, we have an example of an affine semigroup satisfying the conditions of Corollary~\ref{cor:principal}. 

\begin{example}
	Let $\a_1=(1,5), \a_2=(5,1) , \a_3=(2,2), \a_4=(3,3)$. Then
	\begin{eqnarray*}
		\AP(S,E)&=&\{0,\a_3, \a_4, 2\a_3, \a_3+\a_4, 2\a_3+\a_4\}\\
		&=&\{0,\w_1=(2,2),\w_2=(3,3),\w_3=(4,4),\w_4=(5,5),\w_5=(7,7)\},
	\end{eqnarray*}
	and  $r(\Ap(S,E)=\{0,\b_1=(2,2),\b_2=(3,3),\b_3=(4,4),\b_4=(5,5),\b_5=(1,1)\}$.
	As $C_i=\{\w_i\}$, for $i=1,\dots,5$, $\KK[S]$ is Cohen-Macaulay, by Lemma~\ref{CM-Ci}. Moreover, $\max_{\preceq_c}\Ap(S,E)=\{\w_5\}$ and  $\max_{\preceq_c}r(\Ap(S,E))=\{\b_4\}$. Therefore, by Corollary~\ref{cor:principal}, $\fc(S)$ is generated by $\w_5-\b_4=(2,2)$, as an ideal of $\bar{S}=\langle (1,1),(1,5),(5,1)\rangle$. 
\end{example}

For an affine semigroup $S$ with Cohen-Macaulay semigroup ring, if  $\max_{\preceq_c}\Ap(S,E)$ is a singleton, it does not imply that  $\max_{\preceq_c}r(\Ap(S,E))$ is also a singleton.

\begin{example}\label{MR-ex2}
	Let $\a_1=(2,1), \a_2=(1,5), \a_3=(1,1), \a_4=(4,5)$. As the computation in \cite[Example~2]{Mahdavi-Rahmati-2016} shows,   $\max_{\preceq_c}\Ap(S,E)=(6,7)$,  $\max_{\preceq_{c}}r(\Ap(S,E))=\{(2,2),(2,3),(2,4),(2,5)\}$ and $\KK[S]$ is Cohen-Macaulay. Therefore, $\fc(S)$ is generated by $\{(4,5),(4,4),(4,3),(4,2)\}$ as an ideal of $\bar{S}=\langle (1,1),(1,2),(1,3),(1,4),(1,5),$ $(2,1)\rangle$, by Corollary~\ref{cor:principal}.
\end{example}

 Note that  $\max_{\preceq_c}\Ap(S,E)\subseteq\max_{\preceq_S}\Ap(S,E)$. In particular, if $\KK[S]$ is Gorenstein, then $\max_{\preceq_c}\Ap(S,E)$ has a single element. The converse is not true, for instance $\max_{\preceq_S}\Ap(S,E)=\{(6,7),(5,5)\}$ in Example~\ref{MR-ex2}, while $\max_{\preceq_c}\Ap(S,E)$ has a single element. 

\begin{corollary}
			If $\KK[S]$  is a Gorenstein ring and $\max_{\preceq_c}(r(\Ap(S,E)))$ has  a single element, then $\fc(S)$ is a principal ideal of $\bar{S}$ .

\end{corollary}
\begin{remark}
	Let $S\subseteq\NN$ be a numerical semigroup with multiplicity $e$, that is $e=\min(S\setminus\{0\})$. As we mentioned in Example~\ref{F+1},  $\fc(S)$ is generated by $F+1$ as an ideal of $\bar{S}=\NN$, where $F=\max(\ZZ\setminus S)$. Note that $r(\Ap(S,e))=\{0,1,\dots,e-1\}$. As an immediate consequence of Corollary~\ref{cor:principal}, we derive that $F=w-e$, where $w$ is the  maximal number in  $\Ap(S,e)$. This is a fact already proved differently in \cite{selmer}, see also \cite[Theorem~2.12]{RG}.  As Example~\ref{New-Ex} shows,   the conductor of a Cohen-Macaulay affine semigroup $S$ is not necessarily a principal ideal of $\bar{S}$.
\end{remark}

The following is an example of a Cohen-Macaulay simplicial affine semigroup, for which $\max_{\preceq_c}\Ap(S,E)$ is a singleton but  $\fc(S)$ is not principal.  

\begin{example}
	Let $\a_1=(3,0),\a_2=(0,3),\a_3=(2,1)$. Then 	
	 $\AP(S,E)=\{0,\w_1=(2,1),\w_2=(4,2)\}$ and $r(\Ap(S,E))=\{0, \b_1=(2,1), \b_2=(1,2)\}$. 
 Since $C_i=\{\w_i\}$ for $i=1,2$,  $\KK[S]$ is Cohen-Macaulay   by Lemma~\ref{CM-Ci}. Moreover,  $\max_{\preceq_c}\{\w_1,\w_2\}=\{\w_2\}$ and $\max_{\preceq_c}r(\Ap(S,E))=~\{\b_1,\b_2\}$.
By Proposition~\ref{prop:sing-max}, $\fc(S)$ is generated by $\{\w_2-\b_1,\w_2-\b_2\}=\{(2,1),(3,0)\}$ as an ideal of $\bar{S}=\langle (3,0),(0,3),(1,2),(2,1)\rangle$.  
\end{example}

 	\begin{example}\label{4.23}
Let $S$ be the affine semigroup presented in Example~\ref{2.2}. Based on a computation by Macaulay~2 \cite{Mac2}, $\KK[S]$ is Cohen-Macaulay. Note that $\Ap(S,E)$ has two maximal elements $18\a_4+2\a_5=(40,22,20)$ and $16\a_4+4\a_5=(40,24,20)$. Computing the coordinates of these two elements, we find that $\f_{(\w_1,\dots,\w_{k})}=3\a_1+\a_2+3\a_3=(40,23,20)$. One can use GAP \cite{GAP} with an implementation based on Theorem~\ref{fff} (provided to us by the referee), to find that  $\fc(S)$ is minimally generated by 
$\{(29, 15, 14)$, $(29, 16, 14),$ $(29, 18, 15), (30, 15, 14),$ $(30, 16, 14),$ $(31, 15, 15),$ $(31, 16, 14),$ $(31, 18, 16),$ $(33, 18, 17),$ $(33, 19, 15),$ $(35, 18, 18)\}$. 
 		
 	\end{example}
 		
 We close this section by an  example in dimension three.
 
 \begin{example}
 	Let $\a_1=(1,2,1), \a_2=(2,3,1), \a_3=(2,1,3), \a_4=(2,3,2), \a_5=(2,2,2), \a_6=(3,3,3)$. Then $\Ap(S,E)=\{0, \w_1=\a_4, \w_2=\a_5, \w_3=\a_6, \w_4=\a_5+\a_6\}$ and $r(\Ap(S,E))=\{0, \b_1=(1,1,1), \b_2=(2,2,2), \b_3=(3,3,3)\}$.  Moreover, $C_1=\{\w_1=\a_1+\b_1,\w_4=\a_2+\a_3+\b_1=2\a_5+\b_1\}$, $C_2=\{\w_2=\b_2\}$   and $C_3=\{\w_3=\b_3\}$.
 	As $\b_2,\b_3\in S$, we have $\fc(S)=\{\a\in S \ ; \ \a+\b_1\in S\}.$ 	
 	Since $\a_i+\b_1\in S$ for $i\in\{1,4,5,6\}$, we have  $\a_1,\a_4=\a_1+\b_1,\a_5,\a_6=\a_5+\b_1$ are in $\fc(S)$.
 	 	If $\c$ is a minimal generator of $\fc(S)$ that  is not in $\{\a_1,\a_5\}$, then $\c=t\a_2+s\a_3$ for some $t,s\in\NN$. Thus,  $r(\c+\b_1)=\b_1$, and we get  $\c+\b_1=\w+\sum^3_{i=1}l_i\a_i$ for some $\w\in C_1$ and $l_1,l_2,l_3\in\NN$. 
 	 	
 	 	If $\w=\w_1$, then $\c+\b_1=\a_1+\b_1+\sum^3_{i=1}l_i\a_i$, consequently $\c\in\a_1+\bar{S}$, which is a contradiction. 
 	 	
 	 	If $\w=\w_4$, then $\c+\b_1=2\a_5+\b_1+\sum^3_{i=1}l_i\a_i$, which implies  $\c\in\a_5+\bar{S}$,  a contradiction.   	
 	
 	Thus, as an ideal of $\bar{S}$, $\fc(S)$ is minimally generated by $\{\a_1,\a_5\}$.  	   	
 	 	 Note that $\f_{(\w_1,\w_2,\w_3)}=\a_1, \f_{(\w_4,\w_2,\w_3)}=\a_2+\a_3$ and  $\a_5=\a_2+\a_3-\b_2$.  In particular, $\fc(S)$ is generated by $\{\f_{(\w_1,\w_2,\w_3)}, \f_{(\w_4,\w_2,\w_3)}-\b_2\}$ as an ideal of $\bar{S}=\langle \a_1,\a_2,\a_3,\b_1\rangle$. 
 \end{example}

\subsection*{Acknowledgments}

The authors are deeply grateful to the anonymous referee for very careful reading of the manuscript, providing valuable  comments, catching several mistakes and suggesting interesting examples. The referee has also provided the computations we mentioned in Example~\ref{4.23}. The quality of the present paper is very much indebted to the time and expertise devoted by the referee and the  editor.  We also would like to thank Ignacio Ojeda for telling us about the results in \cite{Campillo-Gimenez-2000}. Finally, we acknowledge the use of Macaulay2~(\cite{Mac2}) and  GAP (\cite{GAP}) software for the computations.


\end{document}